\documentclass{amsart}
\usepackage{amssymb,amsmath,amsfonts,amsthm,graphicx,enumerate,amscd,latexsym,curves,multibox}
\usepackage[mathscr]{eucal}
\usepackage{epsfig,epsf,xypic,epic}

\usepackage[T1]{fontenc}
\usepackage[sc]{mathpazo}

\newtheorem{thm}{Theorem}[section]

\newtheorem{prop}{Proposition}[section]
\newtheorem{defn}{Definition}[section]
\newtheorem{nb}{Remark}[section]
\newtheorem{conj}{Conjecture}[section]

\numberwithin{equation}{section}

\begin{document}

\title[SYZ transformations]{On SYZ mirror transformations}
\author[K.-W. Chan]{Kwokwai Chan}
\address{Department of Mathematics, Harvard University,
1 Oxford Street, Cambridge, MA 02138, USA}
\email{kwchan@math.harvard.edu}
\author[N.-C. Leung]{Naichung Conan Leung}
\address{The Institute of Mathematical Sciences and Department of Mathematics,
The Chinese University of Hong Kong, Shatin, Hong Kong}
\email{leung@ims.cuhk.edu.hk}

\begin{abstract}
In this expository paper, we discuss how Fourier-Mukai-type
transformations, which we call SYZ mirror transformations, can be
applied to provide a geometric understanding of the mirror symmetry
phenomena for semi-flat Calabi-Yau manifolds and toric Fano
manifolds. We also speculate the possible applications of these
transformations to other more general settings.
\end{abstract}

\maketitle

\tableofcontents

\section{Introduction}

In 1996, Strominger, Yau and Zaslow suggested, in their
ground-breaking work \cite{SYZ96}, a geometric approach to the
mirror symmetry for Calabi-Yau manifolds. Roughly speaking, the
\emph{Strominger-Yau-Zaslow (SYZ) Conjecture} asserts that any
Calabi-Yau manifold $X$ should admit a fibration by special
Lagrangian tori and the mirror of $X$, which is another Calabi-Yau
manifold $Y$, can be obtained by T-duality, i.e. dualizing the
special Lagrangian torus fibration of $X$. Moreover, the symplectic
geometry (A-model) of $X$ should be interchanged with the complex
geometry (B-model) of $Y$, and vice versa, through fiberwise
Fourier-Mukai-type transformations, suitably modified by quantum
corrections. These transformations are called \emph{SYZ mirror
transformations} and they will be the theme in this article.

Much work has been done on the SYZ Conjecture. Following the work of
Hitchin \cite{Hitchin97}, Leung-Yau-Zaslow \cite{LYZ00} and Leung
\cite{Leung00} explained successfully and neatly the mirror symmetry
for semi-flat Calabi-Yau manifolds by using \emph{semi-flat SYZ
mirror transformations}. These are honest fiberwise real
Fourier-Mukai transformations. The advantage in this case is the
absence of quantum corrections by holomorphic curves and discs. This
is due to the fact that the special Lagrangian torus fibrations on
semi-flat Calabi-Yau manifolds do not admit singularities, and,
accordingly, the bases are smooth affine manifolds.

To deal with general compact Calabi-Yau manifolds, however, one
cannot avoid singularities in Lagrangian torus fibrations, and hence
singularities in the base affine manifolds. Consequently, quantum
corrections will come into play. This necessitates the study of
moduli spaces of special Lagrangian submanifolds and affine
manifolds with singularities, which makes the subject much more
sophisticated and difficult. Nevertheless, the recent progress made
by Gross and Siebert \cite{GS07}, after earlier works of Fukaya
\cite{Fukaya02} and Kontsevich-Soibelman \cite{KS04}, was
doubtlessly a significant step towards establishing the SYZ
Conjecture for general compact Calabi-Yau manifolds.\footnote{We
should mention that the Gross-Siebert program is expected to work
for non-Calabi-Yau manifolds (e.g. Fano manifolds) as well.}

On the other hand, mirror symmetry phenomena have also been observed
for Fano manifolds (and other classes of manifolds or orbifolds as
well). The mirror of a Fano manifold $\bar{X}$ is predicted by
Physicists to be given by a \emph{Landau-Ginzburg model}, which is a
pair $(Y,W)$, consisting of a non-compact K\"{a}hler manifold $Y$
and a holomorphic function $W:Y\rightarrow\mathbb{C}$ called the
\emph{superpotential}. A very important class of examples is
provided by toric Fano manifolds. In this case, the mirror manifold
$Y$ is biholomorphic to (a bounded domain of) $(\mathbb{C}^*)^n$ and
the superpotential $W$ is a Laurent polynomial which can be written
down explicitly. Ample evidences have been found in this toric Fano
case; in particular, Cho and Oh \cite{Cho-Oh03} proved that the
superpotential can be computed in terms of the counting of Maslov
index two holomorphic discs in $\bar{X}$ with boundary in Lagrangian
torus fibers. In \cite{Auroux07}, Auroux applied the SYZ philosophy
to the study of the mirror symmetry for a compact K\"{a}hler
manifold equipped with an anticanonical divisor. This is a
generalization of the mirror symmetry for Fano manifolds, and,
again, the mirror is given by a Landau-Ginzburg model. Auroux also
made an attempt to compute the superpotential in terms of the
counting of holomorphic discs, and analyzed the resulting
wall-crossing phenomena. In \cite{Chan-Leung08}, we studied the
mirror symmetry for toric Fano manifolds, again through the SYZ
approach, and we constructed and applied \textit{SYZ mirror
transformations for toric Fano manifolds} to explain various
geometric results implied by mirror symmetry.\\

A brief explanation of the results in \cite{Chan-Leung08} is now in
order; for more details, see Section~\ref{sec3}. Let $\bar{X}$ be a
toric Fano manifold, i.e. a smooth projective toric variety such
that the anticanonical line bundle $K_{\bar{X}}$ is ample. Let
$\omega_{\bar{X}}$ be a toric K\"{a}hler structure on $\bar{X}$. The
moment map $\mu_{\bar{X}}:\bar{X}\rightarrow\bar{P}$ of the
Hamiltonian $T^n$-action on $(\bar{X},\omega_{\bar{X}})$ is a
natural Lagrangian torus fibration. Here
$\bar{P}\subset\mathbb{R}^n$ is a polytope defining
$(\bar{X},\omega_{\bar{X}})$. The restriction of the moment map to
the open dense $T^n$-orbit $X\cong(\mathbb{C}^*)^n\subset\bar{X}$ is
a Lagrangian torus bundle $\mu_X=\mu_{\bar{X}}|_X:X\rightarrow P$,
where $P$ denotes the interior of the polytope $\bar{P}$. Our first
result in \cite{Chan-Leung08} showed that the mirror manifold $Y$ is
nothing but the \emph{SYZ mirror manifold} of $X$, i.e. the total
space of the torus bundle dual to $\mu_X:X\rightarrow P$ (see
Proposition~\ref{prop3.1}).\footnote{More precisely, the SYZ mirror
manifold is a bounded domain in the mirror manifold $Y$ predicted by
Physicists.} Furthermore, the semi-flat SYZ transformation
$\mathcal{F}^{\textrm{sf}}$ takes the exponential of ($\sqrt{-1}$
times) the symplectic structure $\omega_X=\omega_{\bar{X}}|_X$ on
$X$ to the holomorphic volume form $\Omega_Y$ on
$Y$.\footnote{Throughout this paper, we assume that the B-field is
zero.} Note that $\Omega_Y$ determines a complex structure on $Y$ by
declaring that a 1-form $\alpha$ is a $(1,0)$-form if and only if
$\alpha\lrcorner\Omega_Y=0$. This part of the mirror symmetry does
not involve quantum corrections.

To get the superpotential $W$, however, we need to take into account
the quantum corrections due to the anticanonical toric divisor
$D_\infty=\bar{X}\setminus X$, which we have ignored above. Before
doing that, we first take a digression to a well-known construction.
For a simply connected symplectic manifold $(M,\omega)$, let
$\mathscr{L}M$ be the free loop space, i.e. the space of smooth maps
$\gamma:S^1\rightarrow M$. The symplectic structure on $M$ induces a
symplectic structure on $\mathscr{L}M$ which will also be denoted by
$\omega$. The action functional defined by
$$H(\gamma):=\frac{1}{2\pi}\int_{D_\gamma}\omega,$$
where $D_\gamma$ is a disk contracting $\gamma$, becomes a
well-defined function on the universal covering
$\widetilde{\mathscr{L}M}$ of the free loop space $\mathscr{L}M$.
The group of deck transformations is $H_2(M,\mathbb{Z})$. It is not
hard to see that $H$ is the moment map for the built-in $S^1$-action
on $\widetilde{\mathscr{L}M}$, and the gradient flow lines of $H$
are (pseudo-)holomorphic cylinders if we fix a compatible (almost)
complex structure on $M$. Tentatively, the quantum cohomology (or
Floer cohomology) is the $S^1$-equivariant Morse-Witten cohomology
of the moment map $H$ on $\widetilde{\mathscr{L}M}$. However, the
fact that $\widetilde{\mathscr{L}M}$ is infinite dimensional poses
severe difficulties in implementing this idea.

One of our discoveries in \cite{Chan-Leung08} was that a finite
dimensional subspace of $\mathscr{L}M$ is enough to capture the
quantum corrections and recover the quantum cohomology, in the case
when $M=\bar{X}$ is a toric Fano manifold. Consider the subspace
$LX$ of $\mathscr{L}\bar{X}$ consisting of those loops which are
\emph{geodesic} in the Lagrangian torus fibers (with respect to the
flat metrics) of the moment map
$\mu_{\bar{X}}:\bar{X}\rightarrow\bar{P}$. We consider the function
$\Psi$ on $LX$ defined by $\Psi(\gamma)=\exp(-H(\gamma))$ if
$\gamma$ bounds a Maslov index two holomorphic disc and
$\Psi(\gamma)=0$ otherwise. The function
$\Psi:LX\rightarrow\mathbb{C}$, as an object in the A-model of
$\bar{X}$, turns out to be the mirror of the superpotential $W$. In
\cite{Chan-Leung08}, we constructed the SYZ mirror transformation
$\mathcal{F}$ for the toric Fano manifold $\bar{X}$, and showed that
the SYZ mirror transformation of $\Psi$ is precisely the B-model
superpotential $W$. Moreover, by incorporating the symplectic
structure $\omega_X$ and the holomorphic volume form $\Omega_Y$, we
proved that
\begin{eqnarray*}
\mathcal{F}(e^{\sqrt{-1}\omega_X+\Psi}) & = & e^W\Omega_Y,\\
\mathcal{F}^{-1}(e^W\Omega_Y) & = & e^{\sqrt{-1}\omega_X+\Psi},
\end{eqnarray*}
where $\mathcal{F}^{-1}$ is the inverse SYZ mirror transformation
(see Theorem~\ref{thm3.1}). Hence, the \textit{corrected} symplectic
structure on $X$ and the complex structure on $(Y,W)$ are
interchanged by the SYZ mirror transformation. On the other hand, we
identified the small quantum cohomology ring $QH^*(\bar{X})$ of
$\bar{X}$ with an algebra of functions on $LX$, and realized the
quantum product as a \emph{convolution product} (see
Proposition~\ref{prop3.2}). Then, we showed that the SYZ mirror
transformation $\mathcal{F}$ exhibits a natural isomorphism between
$QH^*(\bar{X})$ and the Jacobian ring $Jac(W)$ of the superpotential
$W$, which \emph{takes the quantum product (now as a convolution
product) to the ordinary product of Laurent polynomials, just as
what classical Fourier series do} (see Theorem~\ref{thm3.2}). We
conclude that the mirror symmetry for toric Fano manifolds is
nothing but a Fourier transformation!\\

The main goal of this article is to popularize the use of SYZ mirror
transformations in exploring mirror symmetry phenomena. In Section
\ref{sec2}, we review the use of semi-flat SYZ mirror
transformations in the study of the mirror symmetry for semi-flat
Calabi-Yau manifolds, where quantum corrections are absent. This is
the toy case which lays the basis for subsequent development in the
investigation of the SYZ Conjecture. Section \ref{sec3} discusses
the mirror symmetry for toric Fano manifolds, where quantum
corrections arise due to the anticanonical toric divisor. Following
\cite{Chan-Leung08}, we demonstrate how to construct and apply SYZ
mirror transformations in this case. The final section contains a
brief discussion of possible generalizations.\\

\noindent\textbf{Acknowledgments.} The authors are grateful to the
organizers of the conference "New developments in Algebraic
Geometry, Integrable Systems and Mirror Symmetry" held in Kyoto
University in January 2008 for giving them an opportunity to
participate in such a stimulating and fruitful event. Thanks are
also due to Hiroshi Iritani and Cheol-Hyun Cho for many useful
discussions. Finally, we thank the referee for several helpful
comments. K.-W. C. was partially supported by Harvard University and
the Croucher Foundation Fellowship. N.-C. L. was partially supported
by RGC grants from the Hong Kong Government.

\section{SYZ mirror transformations without corrections}\label{sec2}

In this section, we review the construction of SYZ mirror
transformations for semi-flat Calabi-Yau manifolds and see how they
were applied in the study of semi-flat mirror symmetry.

\subsection{Semi-flat SYZ mirror transformations}

Denote by $N\cong\mathbb{Z}^n$ a rank-$n$ lattice and
$M=\mbox{Hom}(N,\mathbb{Z})$ the dual lattice. Let $D\subset
M_\mathbb{R}=M\otimes_\mathbb{Z}\mathbb{R}$ be a convex domain.
\footnote{More generally, instead of a convex domain, one may
consider a smooth affine manifold.} Then the tangent bundle
$TD=D\times\sqrt{-1}M_\mathbb{R}$ is naturally a complex manifold
with complex coordinates $x_j+\sqrt{-1}y_j$, $j=1,\ldots,n$, where
$x_1,\ldots,x_n\in\mathbb{R}$ and $y_1,\ldots,y_n\in\mathbb{R}$ are
respectively the base coordinates on $D$ and fiber coordinates on
$M_\mathbb{R}$. We have the standard holomorphic volume form
$\Omega_{TD}=d(x_1+\sqrt{-1}y_1)\wedge\ldots\wedge
d(x_n+\sqrt{-1}y_n)$ on $TD$. By taking fiberwise quotient by the
lattice $M\subset M_\mathbb{R}$, we can compactify the fiber
directions to give the complex manifold
$$Y=TD/M=D\times\sqrt{-1}T_M,$$
where $T_M$ denotes the torus $M_\mathbb{R}/M$. The complex
coordinates on $Y$ are naturally given by
$z_j=\exp(x_j+\sqrt{-1}y_j)$, $j=1,\ldots,n$, where
$y_1,\ldots,y_n\in\mathbb{R}/2\pi\mathbb{Z}$ are now coordinates on
$T_M$. Note that $Y$ is biholomorphic to an open part of
$(\mathbb{C}^*)^n=TM_\mathbb{R}/M$. The projection to $D$ is a torus
bundle, which we denote by $\nu_Y:Y\rightarrow D$. The holomorphic
$n$-form $\Omega_{TD}$ descends to give the holomorphic volume form
$$\Omega_Y=\frac{dz_1}{z_1}\wedge\ldots\wedge\frac{dz_n}{z_n}$$
on $Y$. As mentioned in the introduction, $\Omega_Y$ in turn
determines the complex structure on $Y$: a 1-form $\alpha$ is of
$(1,0)$-type if and only if $\alpha\lrcorner\Omega_Y=0$. Further, if
$\phi$ is an elliptic solution of the \emph{real Monge-Amp\`{e}re
equation}
$$\textrm{det}\Big(\frac{\partial^2\phi}{\partial x_j\partial x_k}\Big)=const,$$
then the K\"{a}hler form
$$\omega_Y:=\sqrt{-1}\partial\bar{\partial}\phi=\sum_{j,k}\phi_{jk}
dx_j\wedge dy_k,$$ with $\phi_{jk}$ denoting
$\frac{\partial^2\phi}{\partial x_j\partial x_k}$, gives a
Calabi-Yau metric on $Y$, and
$$\nu_Y:Y\rightarrow D$$
becomes a special Lagrangian torus bundle (\emph{SYZ fibration}). In
summary, we have the following structures on the complex
$n$-dimensional semi-flat Calabi-Yau manifold $Y$:

\begin{center}
\begin{tabular}{|l|l|} \hline Riemannian metric & $g_Y=
\sum_{j,k}\phi_{jk}(dx_j\otimes dx_k+dy_j\otimes dy_k)$\\
\hline Holomorphic volume form & $\Omega_Y=\bigwedge_{j=1}^n(dx_j+\sqrt{-1}dy_j)$\\
\hline Symplectic form & $\omega_Y=\sum_{j,k}\phi_{jk}dx_j\wedge
dy_k$\\
\hline SYZ fibration & $\nu_Y:Y\rightarrow D$\\
\hline
\end{tabular}
\vspace{0.3cm}
\end{center}

As suggested in the monumental work
Strominger-Yau-Zaslow~\cite{SYZ96}, the mirror of $Y$, which is
another Calabi-Yau manifold we denote by $X$, should be given by the
moduli space of pairs $(L,\nabla)$, where $L$ is a special
Lagrangian torus fiber in $Y$, and $\nabla$ is a flat
$U(1)$-connection on the trivial complex line bundle
$L\times\mathbb{C}\rightarrow L$. This is nothing but the total
space of the torus fibration $\mu_X:X=D\times\sqrt{-1}T_N\rightarrow
D$, where $T_N=N_\mathbb{R}/N=(T_M)^\vee$ and
$N_\mathbb{R}=N\otimes_\mathbb{Z}\mathbb{R}$, which is dual to
$\nu_Y:Y\rightarrow D$. This is called \emph{T-duality} in physics.
Furthermore, $X$ can naturally be viewed as the fiberwise quotient
of the cotangent bundle $T^*D=D\times\sqrt{-1}N_\mathbb{R}$ by the
lattice $N\subset N_\mathbb{R}$. In particular, the standard
symplectic form $\omega_{T^*D}=\sum_{j=1}^n dx_j\wedge du_j$
descends to give a symplectic form
$$\omega_X=\sum_{j=1}^n dx_j\wedge du_j$$
on $X=T^*D/N$, where $u_1,\ldots,u_n\in\mathbb{R}/2\pi\mathbb{Z}$
are coordinates on $T_N$. Through the metric
$$g_X=\sum_{j,k}(\phi_{jk}dx_j\otimes dx_k+\phi^{jk}du_j\otimes du_k),$$
where $(\phi^{jk})$ is the inverse matrix of $(\phi_{jk})$, we
obtain a complex structure on $X$ with complex coordinates given by
$d\log(w_j)=\sum_{k=1}^n\phi_{jk}dx_k+\sqrt{-1}du_j$. There is a
corresponding holomorphic volume form which can be written as
$$\Omega_X=\frac{dw_1}{w_1}\wedge\ldots\wedge\frac{dw_n}{w_n}
=\bigwedge_{j=1}^n(\sum_{k=1}^n\phi_{jk}dx_k+\sqrt{-1}du_j).$$ The
projection map
$$\mu_X:X\rightarrow D$$
now naturally becomes a special Lagrangian torus fibration. In
summary, we have the following structures on $X$:

\begin{center}
\begin{tabular}{|l|l|} \hline Riemannian metric & $g_X=
\sum_{j,k}(\phi_{jk}dx_j\otimes dx_k+\phi^{jk}du_j\otimes du_k)$\\
\hline Holomorphic volume form & $\Omega_X=
\bigwedge_{j=1}^n(\sum_{k=1}^n\phi_{jk}dx_k+\sqrt{-1}du_j)$\\
\hline Symplectic form & $\omega_X=\sum_{j=1}^n dx_j\wedge du_j$\\
\hline SYZ fibration & $\mu_X:X\rightarrow D$\\
\hline
\end{tabular}
\vspace{0.3cm}
\end{center}

We remark that both $Y$ and $X$ admit natural Hamiltonian
$T^n$-actions, but while $\mu:X\rightarrow D$ is a moment map for
the $T_N$-action on $X$, $\nu:Y\rightarrow D$ is not a moment map
for the $T_M$-action on $Y$. In fact, a moment map
$\mu_Y:Y\rightarrow N_\mathbb{R}$ for the $T_M$-action on $Y$ is
given by
$$\mu_Y=L_\phi\circ\nu_Y,$$
where $L_\phi:D\rightarrow N_\mathbb{R}$ is the \emph{Legendre
transform} of $\phi$ defined by
$$L_\phi(x_1,\ldots,x_n)=d\phi_x=\Big(\frac{\partial\phi}{\partial x_1},\ldots,\frac{\partial\phi}{\partial x_n}\Big).$$
Since $\phi$ is convex, the image $D^*=L_\phi(D)$ is an open convex
subset of $(M_\mathbb{R})^*=N_\mathbb{R}$. (For this and other
properties of the Legendre transform, see the book of Guillemin
\cite{Guillemin94}, Appendix 1.) In the action coordinates
$x^1,\ldots,x^n$ of $D^*$, which are given by $\frac{\partial
x^j}{\partial x_k}=\phi_{jk}$, the various structures on $Y$ can be
rewritten as:

\begin{center}
\begin{tabular}{|l|l|} \hline Riemannian metric & $g_Y=
\sum_{j,k}(\phi^{jk}dx^j\otimes dx^k+\phi_{jk}dy_j\otimes dy_k)$\\
\hline Holomorphic volume form & $\Omega_Y=
\bigwedge_{j=1}^n(\sum_{k=1}^n\phi^{jk}dx^k+\sqrt{-1}dy_j)$\\
\hline Symplectic form & $\omega_Y=\sum_{j=1}^n dx^j\wedge dy_j$\\
\hline SYZ fibration & $\mu_Y:Y\rightarrow D^*$\\
\hline
\end{tabular}
\vspace{0.3cm}
\end{center}

We call $X$ the \emph{SYZ mirror manifold} of $Y$ (and vice versa)
since the symplectic (resp. complex) geometry of $X$ and the complex
(resp. symplectic) geometry of $Y$ are interchanged under the
\emph{semi-flat SYZ mirror transformation}, which is described as
follows.

First recall that the dual torus $T_M=(T_N)^\vee$ can be interpreted
as the moduli space of flat $U(1)$-connections on the trivial
complex line bundle over $T_N$. More precisely, given
$y=(y_1,\ldots,y_n)\in M_\mathbb{R}\cong\mathbb{R}^n$, we have a
flat $U(1)$-connection
$$\nabla_y=d+\frac{\sqrt{-1}}{2}\sum_{j=1}^n y_j du_j$$
on $T_N\times\mathbb{C}\rightarrow\mathbb{C}$. The \textit{holonomy}
of $\nabla_y$ is given by the map
$$\textrm{hol}_{\nabla_y}:N\rightarrow U(1),\ v\mapsto
e^{-\sqrt{-1}\langle y,v\rangle}.$$ Hence, $\nabla_y$ is gauge
equivalent to the trivial connection if and only if $y\in
M\cong(2\pi\mathbb{Z})^n$. Moreover this construction gives all flat
$U(1)$-connections on the trivial complex line bundle over $T_N$ up
to unitary gauge transformations. The universal $U(1)$-bundle, i.e.
the Poincar\'{e} line bundle $\mathcal{P}$, is given by the trivial
complex line bundle $(T_N\times T_M)\times\mathbb{C}\rightarrow
T_N\times T_M$ equipped with the connection
$d+\frac{\sqrt{-1}}{2}\sum_{j=1}^n (y_jdu_j-u_jdy_j)$. The curvature
of this connection is the two form
$$F=\sqrt{-1}\sum_{j=1}^n dy_j\wedge du_j.$$

Now consider the relative version of this picture. Let $X\times_D
Y=D\times\sqrt{-1}(T_N\times T_M)$ be the fiber product of the dual
torus bundles $\mu:X\rightarrow D$ and $\nu:Y\rightarrow D$. By
abuse of notations, we still use $\mathcal{P}$ and
$F=\sqrt{-1}\sum_{j=1}^n dy_j\wedge du_j\in\Omega^2(X\times_D Y)$ to
denote the fiberwise universal line bundle and curvature two form
respectively.
\begin{defn}
The semi-flat SYZ mirror transformation
$$\mathcal{F}^{\textrm{sf}}:\Omega^*(X)\rightarrow\Omega^*(Y)$$ is defined by
\begin{eqnarray*}
\mathcal{F}^{\textrm{sf}}(\alpha) & = &
\frac{1}{(2\pi\sqrt{-1})^n}\pi_{Y,*}(\pi_X^*(\alpha)\wedge
e^{\sqrt{-1}F})\\
& = & \frac{1}{(2\pi\sqrt{-1})^n}\int_{T_N}\pi_X^*(\alpha)\wedge
e^{\sqrt{-1}F},
\end{eqnarray*}
where $\pi_X:X\times_D Y\rightarrow X$ and $\pi_Y:X\times_D
Y\rightarrow Y$ are the two projections.
\end{defn}
What is crucial is that this Fourier-Mukai-type transformation
transforms the symplectic structure on $X$ to the complex structure
on $Y$ in the sense of the following two propositions. These already
appeared in \cite{Chan-Leung08}, Proposition 3.2. We include their
proofs, which are somewhat interesting, here for completeness.
\begin{prop}
$$\mathcal{F}^{\textrm{sf}}(e^{\sqrt{-1}\omega_X})=\Omega_Y.$$
\end{prop}
\begin{proof}
\begin{eqnarray*}
\mathcal{F}^{\textrm{sf}}(e^{\sqrt{-1}\omega_X}) & = &
\frac{1}{(2\pi\sqrt{-1})^n}\int_{T_N}\pi_X^*(e^{\sqrt{-1}\omega_X})\wedge e^{\sqrt{-1}F}\\
& = & \frac{1}{(2\pi\sqrt{-1})^n}\int_{T_N} e^{\sqrt{-1}\sum_{j=1}^n (dx_j+\sqrt{-1}dy_j)\wedge du_j}\\
& = & \frac{1}{(2\pi\sqrt{-1})^n}\int_{T_N}\bigwedge_{j=1}^n\big(1+\sqrt{-1}(dx_j+\sqrt{-1}dy_j)\wedge du_j\big)\\
& = & \frac{1}{(2\pi)^n}\int_{T_N}
\Bigg(\bigwedge_{j=1}^n(dx_j+\sqrt{-1}dy_j)\Bigg)\wedge
du_1\wedge\ldots\wedge du_n\\
& = & \Omega_Y,
\end{eqnarray*}
where we have $\int_{T_N}du_1\wedge\ldots\wedge du_n=(2\pi)^n$ in
the final step.
\end{proof}
\noindent As a mirror transformation, $\mathcal{F}^{\textrm{sf}}$
should have the \emph{inversion property}. This is the following
proposition.
\begin{prop}
If we define the inverse transform
$(\mathcal{F}^{\textrm{sf}})^{-1}:\Omega^*(Y)\rightarrow\Omega^*(X)$
by
\begin{eqnarray*}
(\mathcal{F}^{\textrm{sf}})^{-1}(\alpha) & = &
\frac{1}{(2\pi\sqrt{-1})^n}\pi_{X,*}(\pi_Y^*(\alpha)\wedge
e^{-\sqrt{-1}F})\\
& = & \frac{1}{(2\pi\sqrt{-1})^n}\int_{T_M}\pi_Y^*(\alpha)\wedge
e^{-\sqrt{-1}F},
\end{eqnarray*}
then we have
$$(\mathcal{F}^{\textrm{sf}})^{-1}(\Omega_Y)=e^{\sqrt{-1}\omega_X}.$$
\end{prop}
\begin{proof}
\begin{eqnarray*}
(\mathcal{F}^{\textrm{sf}})^{-1}(\Omega_Y) & = &
\frac{1}{(2\pi\sqrt{-1})^n}\int_{T_M}\pi_Y^*(\Omega_Y)\wedge e^{-\sqrt{-1}F}\\
& = & \frac{1}{(2\pi\sqrt{-1})^n}\int_{T_M}
\Bigg(\bigwedge_{j=1}^n(dx_j+\sqrt{-1}dy_j)\Bigg)\wedge
e^{\sum_{j=1}^n dy_j\wedge du_j}\\
& = & \frac{1}{(2\pi\sqrt{-1})^n}\int_{T_M}
\bigwedge_{j=1}^n\big((dx_j+\sqrt{-1}dy_j)
\wedge e^{dy_j\wedge du_j}\big)\\
& = & \frac{1}{(2\pi\sqrt{-1})^n}\int_{T_M}\bigwedge_{j=1}^n
\big(dx_j+\sqrt{-1}dy_j+dx_j\wedge dy_j\wedge du_j\big)\\
& = &
\frac{1}{(2\pi)^n}\int_{T_M}\bigwedge_{j=1}^n\big((1+\sqrt{-1}dx_j\wedge
du_j)\wedge dy_j\big)\\
& = & \frac{1}{(2\pi)^n}\int_{T_M}\bigwedge_{j=1}^n\big(e^{\sqrt{-1}dx_j\wedge du_j}\wedge dy_j\big)\\
& = & \frac{1}{(2\pi)^n}\int_{T_M}e^{\sqrt{-1}\sum_{j=1}^n dx_j\wedge du_j}\wedge dy_1\wedge\ldots\wedge dy_n\\
& = & e^{\sqrt{-1}\omega_X}.
\end{eqnarray*}
\end{proof}
By exactly the same arguments, one can also show that
$$\mathcal{F}^{\textrm{sf}}(\Omega_X)=e^{\sqrt{-1}\omega_Y},
\
(\mathcal{F}^{\textrm{sf}})^{-1}(e^{\sqrt{-1}\omega_Y})=\Omega_X.$$
If we take into account the \textit{B-fields}, then the semi-flat
SYZ transformation will give an identification between the moduli
space of complexified K\"{a}hler structures on $X$ with the moduli
space of complex structures on $Y$, and vice versa. For this and
transformations of other geometric structures, we refer the reader
to Leung \cite{Leung00}.

\subsection{Transformations of branes}

Lying at the heart of the SYZ Conjecture is the basic but important
observation that a point $z=\exp(x+\sqrt{-1}y)\in Y$ defines a flat
$U(1)$-connection $\nabla_y$ on the trivial complex line bundle over
the special Lagrangian torus fiber $L_x=\mu_X^{-1}(x)$. Now, the
point $z\in Y$ together with its structure sheaf $\mathcal{O}_z$ can
be considered as a \emph{B-brane} on $Y$; while the pair
$(L_x,\mathbb{L}_y)$, where $\mathbb{L}_y$ denotes the flat
$U(1)$-bundle $(L_x\times\mathbb{C},\nabla_y)$, gives an
\emph{A-brane} on $X$. This implements the simplest case of
correspondence between branes on mirror manifolds via SYZ
transformations:
$$(L_x,\mathbb{L}_y)\longleftrightarrow(z,\mathcal{O}_z).$$
The space of infinitestimal deformations of the A-brane
$(L_x,\mathbb{L}_y)$, which is given by $H^1(L_x,\mathbb{R})\times
H^1(L_x,\sqrt{-1}\mathbb{R})=H^1(L_x,\mathbb{C})$, is canonically
identified with the tangent space $T_zY$, the space of
infinitestimal deformations of the sheaf $\mathcal{O}_z$.\\

On the other hand, consider a section $L=\{(x,u(x))\in X:x\in D\}$
of $\mu_X:X\rightarrow D$. The submanifold $L$ is Lagrangian if and
only if (locally) there exists a function $f$ such that
$u_j=\frac{\partial f}{\partial x_j}$. By the above observation (now
used in the opposite way), a point $(x,u(x))\in L$ determines a flat
$U(1)$-connection $\nabla_{u(x)}$ on the trivial complex line bundle
over the fiber $(L_x)^\vee=\nu_Y^{-1}(x)$. The family of points
$\{(x,u(x)):x\in D\}$ thus patch together to give the
$U(1)$-connection
$$\nabla_L=d_Y-\frac{\sqrt{-1}}{2}\sum_{j=1}^nu_j(x)dy_j$$
on a certain complex line bundle over $Y$; its curvature two form is
given by
$$F_L=d_Y\Big(-\frac{\sqrt{-1}}{2}\sum_{j=1}^nu_j(x)dy_j\Big)
=-\frac{\sqrt{-1}}{2}\sum_{j,k}\frac{\partial u_j}{\partial
x_k}dx_k\wedge dy_j,$$ and, in particular,
$$F_L^{2,0}=\frac{1}{8}\sum_{j<k}\Big(\frac{\partial u_j}{\partial x_k}
-\frac{\partial u_k}{\partial x_j}\Big)\frac{dz_j}{z_j}\wedge
\frac{dz_k}{z_k}.$$ We conclude that $\nabla_L$ is integrable, i.e.
$F_L^{2,0}=0$, if and only if $L$ is Lagrangian. More generally, we
can equip $L$ with a flat $U(1)$-bundle
$\mathbb{L}=(L\times\mathbb{C},d_L+\alpha)$, where
$\alpha\in\Omega^1(L,\mathbb{R})$ is a closed (and hence exact)
one-form. The A-brane $(L,\mathbb{L})$ is then transformed to the
$U(1)$-connection
$$\nabla_{L,\mathbb{L}}=\nabla_L+\alpha,$$
which again is integrable if and only if $L$ is Lagrangian.
Furthermore, one can prove that $\nabla_{L,\mathbb{L}}$ satisfies
the \emph{deformed Hermitian-Yang-Mills equations} if and only if
$L$ is special Lagrangian (see Leung-Yau-Zaslow \cite{LYZ00} and
Leung \cite{Leung00} for the detailed proofs).
$\nabla_{L,\mathbb{L}}$ is a connection on the holomorphic line
bundle over $Y$ given by the semi-flat SYZ transformation of
$\mathbb{L}$:
$$\mathcal{L}_{L,\mathbb{L}}=\pi_{Y,*}(\pi_X^*(\iota_*\mathbb{L})\otimes\mathcal{P}),$$
where $\iota:L\hookrightarrow X$ is the inclusion map. In
conclusion, the A-brane $(L,\mathbb{L})$ is transformed to the
B-brane $(Y,\mathcal{L}_{L,\mathbb{L}})$ through semi-flat SYZ
transformations:
$$(L,\mathbb{L})\longleftrightarrow(Y,\mathcal{L}_{L,\mathbb{L}}).$$

\section{SYZ mirror transformations with corrections}\label{sec3}

In the previous section, we see that T-duality and SYZ mirror
transformations can be applied successfully to give a geometric
understanding of the mirror symmetry for semi-flat Calabi-Yau
manifolds. However, no quantum corrections were involved in this
case due to the absence of holomorphic curves and discs. The
existence of quantum corrections is also closely related to the
singularities of the Lagrangian torus fibrations, which again are
not present in the semi-flat case. In this section, following
\cite{Chan-Leung08}, we are going to discuss how SYZ mirror
transformations can be applied to a case where quantum corrections
\emph{do} exist, namely, the mirror symmetry for toric Fano
manifolds.

\subsection{Mirror symmetry for toric Fano manifolds}

We begin with a more detailed description of the mirror picture for
toric Fano manifolds \cite{Givental94}, \cite{Kontsevich98},
\cite{Hori-Vafa00}. Let $\bar{P}\subset M_\mathbb{R}$ be a smooth
reflexive polytope given by the inequalities
$$\langle x,v_i\rangle\geq\lambda_i,\quad i=1,\ldots,d,$$
where $v_1,\ldots,v_d\in N$ are primitive vectors and
$\langle\cdot,\cdot\rangle:M_\mathbb{R}\times
N_\mathbb{R}\rightarrow\mathbb{R}$ is the dual pairing. This
determines a toric Fano manifold $\bar{X}$, together with a
K\"{a}hler structure $\omega_{\bar{X}}$. Unlike the case of
Calabi-Yau manifolds, the mirror of $\bar{X}$ is not another compact
K\"{a}hler manifold, but a Landau-Ginzburg model: a pair $(Y,W)$
consisting of a noncompact K\"{a}hler manifold $Y$, which (as a
complex manifold) is biholomorphic to (a bounded domain of)
$(\mathbb{C}^*)^n$, and the Laurent polynomial
$$W=e^{\lambda_1}z^{v_1}+\ldots+e^{\lambda_d}z^{v_d}:Y\rightarrow\mathbb{C},$$
which is called the superpotential. Here $z^{v_i}$ denotes the
monomial $z_1^{v_{i1}}\ldots z_n^{v_{in}}$ in the coordinates
$z_1,\ldots,z_n$ of $Y$. For example, if
$P=\{(x_1,x_2,x_3)\in\mathbb{R}^3:x_1\geq0, x_2\geq0, x_1+x_2\leq
t\}$, then $\bar{X}=\mathbb{C}P^2$ and the mirror Landau-Ginzburg
model is given by the Laurent polynomial
$W(z_1,z_2)=z_1+z_2+\frac{e^{-t}}{z_1z_2}$ on $Y=(\mathbb{C}^*)^2$.

Among the many mirror symmetry predictions are the following
conjectures:
\begin{conj}$\mbox{}$
\begin{enumerate}
\item[1.] The small quantum cohomology ring $QH^*(\bar{X})$ of $\bar{X}$ is
isomorphic to the Jacobian ring $Jac(W)$ of $W$, where
$$Jac(W)=\mathbb{C}[z_1^{\pm1},\ldots,z_n^{\pm1}]/\langle
\partial_1W,\ldots,\partial_nW\rangle,$$
and $\partial_j$ denotes $z_j\frac{\partial}{\partial z_j}$.

\item[2.] (Homological mirror symmetry, see \cite{Kontsevich98},
\cite{Seidel00}, \cite{Orlov03}) There are equivalences of
triangulated categories
\begin{eqnarray*}
D^bCoh(\bar{X}) & \cong & D^\pi Fuk(Y,W)\\
D^\pi Fuk(\bar{X}) & \cong & D_{Sing}(Y,W)
\end{eqnarray*}
where $D^\pi Fuk(Y,W)$ is (a suitably defined version of) the
derived Fukaya category of the Landau-Ginzurg model $(Y,W)$ and
$D_{Sing}(Y,W)$ is the category of singularities of $(Y,W)$.
\end{enumerate}
\end{conj}
Substantial evidences \cite{Givental97}, \cite{HIV00},
\cite{Seidel00}, \cite{Ueda04}, \cite{AKO04}, \cite{AKO05},
\cite{Abouzaid05}, \cite{Abouzaid06}, \cite{Cho-Oh03}, \cite{Cho04}
have been found for these conjectures, while evidence in the
Calabi-Yau and other non-toric cases is much rarer. This is partly
due to the fact that geometric structures on toric varieties are
highly computable and explicit, making them an exceptionally fertile
testing ground for techniques and conjectures.

One of these explicit structures: the Lagrangian torus fibration on
$\bar{X}$ given by the moment map
$\mu_{\bar{X}}:\bar{X}\rightarrow\bar{P}$ of the Hamiltonian
$T_N$-action on $(\bar{X},\omega_{\bar{X}})$, is particularly
important in the SYZ spproach and in the constructions of SYZ mirror
transformations. Let
$$\mu_X:X\rightarrow P$$
be the restriction of the moment map to the open dense $T_N$-orbit
$X=\bar{X}\setminus D_\infty$, where $D_\infty=\bigcup_{i=1}^d D_i$
is the anticanonical toric divisor, and $P$ is the interior of
$\bar{P}$. In the symplectic (or action-angle) coordinates,
$$X=T^*P/N=P\times\sqrt{-1}T_N$$
and the restriction of $\omega_{\bar{X}}$ to $X$ is nothing but the
standard symplectic structure
$$\omega_X=\sum_{j=1}^n dx_j\wedge du_j,$$
where $x_1,\ldots,x_n\in\mathbb{R}$ and
$u_1,\ldots,u_n\in\mathbb{R}/2\pi\mathbb{Z}$ are respectively the
base coordinates on $P$ and fiber coordinates on $T_N$ (see Abreu
\cite{Abreu00}). Now we are in exactly the same situation as in the
previous section and it is tempting to assert that the mirror
manifold $Y$ predicted by Physicists is given by the SYZ mirror
manifold of $X$, which is $TP/M=P\times\sqrt{-1}T_M$. This is indeed
nearly the case.
\begin{prop}[Proposition 3.1 in \cite{Chan-Leung08}]\label{prop3.1}
The mirror manifold $Y=(\mathbb{C}^*)^n$ predicted by Physicists
contains the SYZ mirror manifold $TP/M=P\times\sqrt{-1}T_M$ of
$X=\bar{X}\setminus D_\infty$ as a bounded domain
$$\{(z_1,\ldots,z_n)\in Y:|e^{\lambda_i}z^{v_i}|<1\textrm{ for }i=1,\ldots,d\}.$$
Equivalently, the SYZ mirror manifold is given by the preimage of
$P\subset M_\mathbb{R}=\mathbb{R}^n$ under the Log map
$$\textrm{Log}:(\mathbb{C}^*)^n\rightarrow\mathbb{R}^n,\
(z_1,\ldots,z_n)\mapsto(\log |z_1|,\ldots,\log |z_n|).$$
\end{prop}

The same result also appeared in Auroux's paper \cite{Auroux07}
(Proposition 4.2). Also included in his paper was a discussion of
the issue that the SYZ mirror manifold (a bounded domain in
$(\mathbb{C}^*)^n$) is "smaller" than Hori-Vafa's mirror manifold
(the whole $(\mathbb{C}^*)^n$). There is evidence (say, in
Abouzaid's works \cite{Abouzaid05}, \cite{Abouzaid06}) showing that
one should work with the SYZ mirror manifold, instead of the whole
$(\mathbb{C}^*)^n$, in studying mirror symmetry. In any case, we
will use and work with the SYZ mirror manifold, i.e. the bounded
domain in $(\mathbb{C}^*)^n$, and denote it by $Y$ henceforth.

In terms of the coordinates
$z_1=\exp(-x_1-\sqrt{-1}y_1),\ldots,z_n=\exp(-x_n-\sqrt{-1}y_n)\in\mathbb{C}^*$
of $Y\subset(\mathbb{C}^*)^n$, the holomorphic volume form is given
by the standard one on $(\mathbb{C}^*)^n$:
$$\Omega_Y=\frac{dz_1}{z_1}\wedge\ldots\wedge\frac{dz_n}{z_n}$$
and the torus fibration $\nu_Y:Y\rightarrow P$ is the restriction of
the Log map. We remark that metrically we are \emph{not} considering
$X=\bar{X}\setminus D_\infty$ as a Calabi-Yau manifold; instead of
the semi-flat Calabi-Yau metric, we use the $T_N$-invariant
K\"{a}hler metric on $\bar{X}$ (and the corresponding dual metric on
$Y$). These are defined (cf. Guillemin \cite{Guillemin94} and Abreu
\cite{Abreu00}) using the strictly convex function
$\phi_P:P\rightarrow\mathbb{R}$ given by
$$\phi_P(x)=\frac{1}{2}\sum_{i=1}^d l_i(x)\log l_i(x),$$
where $l_i(x)=\langle x,v_i\rangle-\lambda_i$ for $i=1,\ldots,d$,
instead of a solution of the real Monge-Amp\`{e}re equation. For
example, this gives the standard Fubini-Study metric on
$\bar{X}=\mathbb{C}P^n$. Using these metrics and the corresponding
holomorphic volume forms, $X$ and $Y$ are almost Calabi-Yau
manifolds and the torus fibers of $\mu_X$ and $\nu_Y$ are special
Lagrangian submanifolds (also see Section 2 in Auroux
\cite{Auroux07}).

\subsection{SYZ transformations for toric Fano manifolds}

By applying the semi-flat SYZ mirror transformation or T-duality, we
can obtain the mirror manifold $Y$. But where comes the
superpotential $W:Y\rightarrow\mathbb{C}$? Recall that, in applying
T-duality, we have completely ignored the compactification of $X$,
which is given by adding the anticanonical toric divisor
$D_\infty=\bigcup_{i=1}^d D_i$. As suggested in the foundational
work of Fukaya-Oh-Ohta-Ono \cite{FOOO00}, this has tremendous effect
on the Floer theory of the Lagrangian torus fibers of
$\mu_X:X\rightarrow P$, and this is indeed where quantum corrections
by holomorphic discs come into play.

As have been discussed in the introduction, motivated by the idea of
using Morse theory on the free loop space $\mathscr{L}\bar{X}$ to
construct the quantum cohomology $QH^*(\bar{X})$, we introduce the
subspace $LX\subset\mathscr{L}\bar{X}$ consisting of those loops
which are geodesic in the Lagrangian torus fibers of the moment map
$\mu_X:X\rightarrow P$, i.e.
$$LX=\{\gamma\in\mathscr{L}\bar{X}:\textrm{$\gamma$ is a geodesic in
$L_x=\mu_X^{-1}(x)$ for some $x\in P$}\}.$$ Concretely, we have
$$LX=X\times N=P\times\sqrt{-1}T_N\times N,$$
and we consider it as a (trivial) $\mathbb{Z}^n$-cover of $X$,
$\pi:LX\rightarrow X$. Notice that, for each Lagrangian torus fiber
$L_x$, $x\in P$, we have a canonical identification $\pi_1(L_x)\cong
N$.

We are going to define a function $\Psi$ on $LX$ in terms of the
counting of holomorphic discs in $\bar{X}$ of minimal Maslov index.
This will recapture the information of the compactification of $X$
by $D_\infty$, which we have ignored previously, and $\Psi$ serves
as the object in the A-model of $\bar{X}$ mirror to the
superpotential $W$. To do this, let's first recall the fundamental
results of Cho-Oh \cite{Cho-Oh03} on the classification of
holomorphic discs in $\bar{X}$ with boundary in Lagrangian torus
fibers of $\mu_X:X\rightarrow P$.

Let $L_x=\mu_X^{-1}(x)$ be the Lagrangian torus fiber in $X$ over a
point $x\in P$. Then the relative homotopy group
$\pi_2(\bar{X},L_x)$ is generated by the Maslox index two classes
$\beta_1,\ldots,\beta_d$, which are represented by holomorphic discs
in $(\bar{X},L_x)$. Note that we have, $\partial\beta_i=v_i$, for
$i=1,\ldots,d$, where
$\partial:\pi_2(\bar{X},L_x)\rightarrow\pi_1(L_x)\cong N$ is the
natural boundary map. In \cite{Cho-Oh03}, Cho and Oh proved that,
for $i=1,\ldots,d$ and for each point $p\in L_x$, there is a unique
(up to automorphism of the domain) Maslov index two $J$-holomorphic
disc $\varphi_i:(D^2,\partial D^2)\rightarrow(\bar{X},L_x)$ in the
class $\beta_i$ which passes through $p$ and intersects the toric
divisor $D_i$ at an interior point.\footnote{Another way to state
this result is the following. Let $\mathscr{M}_1(\beta_i)$ be the
moduli space of $J$-holomorphic discs $\varphi:(D^2,\partial
D^2)\rightarrow(\bar{X},L_x)$ in the class $\beta_i$ with 1 boundary
marked point. Let $ev:\mathscr{M}_1(\beta_i)\rightarrow L_x$ be the
evaluation map at the boundary marked point. Then the result of Cho
and Oh says that $ev_*[\mathscr{M}_1(\beta_i)]=[L_x]$ as $n$-cycles
in $L_x$. See also Sections 3.1 and 4 in Auroux \cite{Auroux07}.}
Here $J$ is the complex structure on $\bar{X}$ determined by the fan
$\Sigma$ dual to $\bar{P}$.
\begin{defn}
For $i=1,\ldots,d$, define $\Psi_i:LX\rightarrow\mathbb{R}$ by
$$\Psi_i(p,v)=\left\{ \begin{array}{ll}
                      n_i(p)\exp(-\frac{1}{2\pi}\int_{\beta_i}\omega_{\bar{X}}) & \textrm{if $v=v_i$}\\
                      0                         & \textrm{if $v\neq v_i$,}
                      \end{array}\right.$$
for $(p,v)\in LX=X\times N$, where $n_i(p)$ is the algebraic number
of Maslov index two $J$-holomorphic discs in
$(\bar{X},L_{\mu_X(p)})$ in the class $\beta_i$ which pass through
$p$. Then set
$$\Psi=\Psi_1+\ldots+\Psi_d:LX\rightarrow\mathbb{R}.$$
\end{defn}
By their definitons, the $T_N$-invariant functions
$\Psi_1,\ldots,\Psi_d$ carry enumerative meaning, although by Cho
and Oh's result, we always have $n_i(p)=1$, for all $i$ and any $p$.
One may think of the $T_N$-invariant function $\Psi$ as recording
which cycle $v\in N=\pi_1(L_x)$ collapses to a point as one goes
towards the anticanonical toric divisor $D_\infty$, or equivalently,
which geodesic loop $\gamma\in LX$ bounds a holomorphic disc of
Maslov index two.
\begin{nb}
Before showing how to transform $\Psi$ to get the superpotential
$W$, we remark that the $T_N$-invariant function
$\Phi:LX\rightarrow\mathbb{R}$ introduced in \cite{Chan-Leung08},
Definition 2.1, is nothing but the "exponential" of $\Psi$, i.e.
$$\Phi=\textrm{Exp }\Psi,$$
where $\textrm{Exp }\Psi$ is defined as
$\sum_{k=0}^\infty\frac{1}{k!}\underbrace{\Psi\star\ldots\star\Psi}_\textrm{$k$
times}$ in which $\star$ denotes the convolution product of a
certain class of functions on $LX$ with respect to the lattice $N$.
Now each point $q=(q_1,\ldots,q_l)$ ($l=d-n$) in the K\"{a}hler cone
$\mathcal{K}(\bar{X})\subset H^2(\bar{X},\mathbb{R})$ determines a
symplectic structure $\omega_{\bar{X}}$ on $\bar{X}$ and we can
choose the polytope $\bar{P}=\{x\in M_\mathbb{R}:\langle
x,v_i\rangle\geq\lambda_i,\quad i=1,\ldots,d\}$ such that
$v_1=e_1,\ldots,v_n=e_n$ is the standard basis of $N=\mathbb{Z}^n$ ,
$\lambda_1=\ldots=\lambda_n=0$ and $\lambda_{n+a}=\log q_a$ for
$a=1,\ldots,l$. We thus get two families of functions
$\{\Psi_q\}_{q\in\mathcal{K}}$ and $\{\Phi_q\}_{q\in\mathcal{K}}$.
By the symplectic area formula of Cho-Oh (\cite{Cho-Oh03}, Theorem
8.1), we have
\begin{equation*}
\int_{D^2}\varphi_i^*\omega_{\bar{X}}=\int_{\beta_i}\omega_{\bar{X}}=2\pi(\langle
x,v_i\rangle-\lambda_i),
\end{equation*}
for $i=1,\ldots,d$. Hence, for any $(p,v)\in LX$,
$$\Psi_i(p,v)=\left\{ \begin{array}{ll}
                      e^{-\langle x,v_i\rangle} & \textrm{if $v=v_i$}\\
                      0                         & \textrm{if $v\neq v_i,$}
                      \end{array}\right.$$
for $i=1,\ldots,n$, and
$$\Psi_{n+a}(p,v)=\left\{ \begin{array}{ll}
                      q_ae^{-\langle x,v_{n+a}\rangle} & \textrm{if $v=v_{n+1}$}\\
                      0                         & \textrm{if $v\neq v_{n+a}$,}
                      \end{array}\right.$$
for $a=1,\ldots,l$, where $x=\mu_X(p)$. It follows that
$$q_a\frac{\partial\Phi_q}{\partial q_a}=\Phi_q\star\Psi_{n+a}$$
for $a=1,\ldots,l$, which is the first part of Proposition 1.1 in
\cite{Chan-Leung08}.
\end{nb}
On the other hand, the functions $\Psi_1,\ldots,\Psi_d$ are
intimately related to the small quantum cohomology $QH^*(\bar{X})$
of $\bar{X}$, as was shown in the following
\begin{prop}[Second part of Proposition
1.1 in \cite{Chan-Leung08}]\label{prop3.2} Assume that $\bar{X}$ is
a product of projective spaces. Then we have a natural isomorphism
of $\mathbb{C}$-algebras
$$QH^*(\bar{X})\cong\mathbb{C}[\Psi_1^{\pm1},\ldots,\Psi_n^{\pm1}]/\mathfrak{L}$$
where $\mathbb{C}[\Psi_1^{\pm1},\ldots,\Psi_n^{\pm1}]$ is the
polynomial algebra generated by $\Psi_1^{\pm1},\ldots,\Psi_n^{\pm1}$
with respect to the convolution product $\star$, and $\mathfrak{L}$
is the ideal generated by linear relations: $\sum_{i=1}^d
a_i\Psi_i\sim\sum_{i=1}^d b_i\Psi_i$ if and only if the
corresponding divisors $\sum_{i=1}^d a_iD_i$ and $\sum_{i=1}^d
b_iD_i$ are linearly equivalent.
\end{prop}
\begin{nb}
By employing Givental's mirror theorem~\cite{Givental97}, one can in
fact show that the proposition holds for \textit{all} toric Fano
manifolds. See Remark 2.3 in \cite{Chan-Leung08} for details.
\end{nb}
We need the assumption that $\bar{X}$ is a product of projective
spaces as we are intended for a geometric understanding of the
isomorphism in Proposition~\ref{prop3.2} by using \textit{tropical
geometry}. This is briefly described as follows (see Subsection 2.2
in \cite{Chan-Leung08} for details). One first defines a tropical
version $QH^*_{trop}(\bar{X})$ of the small quantum cohomology ring
of $\bar{X}$. Since $\bar{X}$ is a product of projective spaces, we
have a one-to-one correspondences between the $J$-holomorphic curves
in $\bar{X}$ which have contribution to the quantum product in
$QH^*(\bar{X})$ and those tropical curves in $N_\mathbb{R}$ which
have contribution to the tropical quantum product in
$QH^*_{trop}(\bar{X})$, by the \textit{correspondence theorem} of
Mikhalkin \cite{Mikhalkin03} and Nishinou-Siebert \cite{NS04}. From
this follows the canonical isomorphism
$$QH^*(\bar{X})\cong QH^*_{trop}(\bar{X}).$$
Then comes a simple but important observation: \textit{Each tropical
curve which has contribution to the tropical quantum product in
$QH^*_{trop}(\bar{X})$ is obtained by gluing tropical discs in
$N_\mathbb{R}$}.\footnote{This idea was recently generalized by
Gross \cite{Gross09} to understand tropically the big quantum
cohomology and mirror symmetry of $\mathbb{C}P^2$.} On the other
hand, these tropical discs are exactly corresponding to the families
of Maslov index two $J$-holomorphic discs in $\bar{X}$ with boundary
in Lagrangian torus fibers, which were used to define the functions
$\Psi_1,\ldots,\Psi_d$. Hence, we naturally have another canonical
isomorphism
$$QH^*_{trop}(\bar{X})\cong\mathbb{C}[\Psi_1^{\pm1},\ldots,\Psi_n^{\pm1}]/\mathfrak{L}.$$
For example, let us take a look at the case of
$\bar{X}=\mathbb{C}P^2$. See Figure 3.1 below.
\begin{figure}[ht]
\setlength{\unitlength}{1mm}
\begin{picture}(100,35)
\put(10,2){\vector(0,1){35}} \put(10,2){\vector(1,0){35}}
\curve(10,32, 40,2) \curve(40,2.8, 40,1.2) \put(40,-1){$t$}
\curve(10.8,32, 9.2,32) \put(7.8,31){$t$} \put(8,0){0}
\put(13,11){$\bar{P}\subset M_\mathbb{R}$} \put(5.3,14.5){$D_1$}
\put(23,-1.3){$D_2$} \put(24.9,17.9){$D_3$} \curve(75,17, 93,17)
\curve(75,17, 75,35) \curve(75,17, 60,2)
\put(74.15,16.15){$\bullet$} \put(71,16){$\xi$} \put(87,18){$v_1$}
\put(75.5,30){$v_2$} \put(67,6.3){$v_3$} \put(75,4){$\Gamma\subset
N_\mathbb{R}$} \put(44,-3){Figure 3.1}
\end{picture}
\end{figure}

Denote by $\{e_1,e_2\}$ the standard basis of $N=\mathbb{Z}^2$. We
have $v_1=(1,0),v_2=(0,1),v_3=(-1,-1)$, and the polytope
$\bar{P}\subset M_\mathbb{R}\cong\mathbb{R}^2$ is defined by the
inequalities
$$x_1\geq0,\ x_2\geq0,\ x_1+x_2\leq t,$$
where $t>0$. There are three toric divisors $D_1,D_2,D_3$
corresponding to three functions $\Psi_1,\Psi_2,\Psi_3\in
C^\infty(LX)$ defined by
\begin{eqnarray*}
\Psi_1(p,v) & = & \left\{ \begin{array}{ll}
                         e^{-x_1} & \textrm{if $v=(1,0)$}\\
                         0 & \textrm{otherwise,}
                         \end{array} \right.\\
\Psi_2(p,v) & = & \left\{ \begin{array}{ll}
                         e^{-x_2} & \textrm{if $v=(0,1)$}\\
                         0 & \textrm{otherwise,}
                         \end{array} \right.\\
\Psi_3(p,v) & = & \left\{ \begin{array}{ll}
                         e^{-(t-x_1-x_2)} & \textrm{if $v=(-1,-1)$}\\
                         0 & \textrm{otherwise,}
                         \end{array} \right.
\end{eqnarray*}
for $(p,v)\in LX$ and $(x_1,x_2)=\mu_X(p)\in P$, respectively. The
small quantum cohomology ring is given by
\begin{eqnarray*}
QH^*(\mathbb{C}P^2) & = & \mathbb{C}[D_1,D_2,D_3]\big/\big\langle
D_1-D_3,D_2-D_3,D_1\ast D_2\ast D_3-q\big\rangle\\
& = & \mathbb{C}[H]\big/\big\langle H^3-q\big\rangle,
\end{eqnarray*}
where we have, by abuse of notations, also use $D_i\in
H^2(\mathbb{C}P^2,\mathbb{C})$ to denote the cohomology class
Poincar\'{e} dual to $D_i$, $H\in H^2(\mathbb{C}P^2,\mathbb{C})$ is
the hyperplane class and $q=e^{-t}$. Fix any point
$p\in\mathbb{C}P^2\setminus D_\infty$, then the quantum corrections,
which appear in the relation
$$D_1\ast D_2\ast D_3=H^3=q,$$
is due to the unique holomorphic curve
$\varphi:(\mathbb{P}^1;x_1,x_2,x_3,x_4)\rightarrow\mathbb{C}P^2$ of
degree 1 (i.e. a line) with 4 marked points such that
$\varphi(x_4)=p$ and $\varphi(x_i)\in D_i$, for $i=1,2,3$. Let
$x=\mu_X(p)\in P$ and $L_x=\mu_X^{-1}(x)$ be the Lagrangian torus
fiber containing $p$. Using tropical geometry, one sees that there
is a tropical curve $\Gamma$ in $N_\mathbb{R}$ with three unbounded
edges in the directions $v_1,v_2,v_3$ and the vertex mapped to
$\xi=\textrm{Log}(p)\in N_\mathbb{R}$, which is corresponding to
this holomorphic curve (see Figure 3.1 above). Here, we identify $X$
with $(\mathbb{C}^*)^2$, and
$\textrm{Log}:X=(\mathbb{C}^*)^2\rightarrow
N_\mathbb{R}=\mathbb{R}^2$ is the Log map we defined in
Proposition~\ref{prop3.1}.
\begin{figure}[ht]
\setlength{\unitlength}{1mm}
\begin{picture}(100,35)
\curve(20,17, 38,17) \curve(20,17, 20,35) \curve(20,17, 5,2)
\put(10,1){$\Gamma$} \put(19.1,16.1){$\bullet$} \put(16,16){$\xi$}
\put(42,16){glued from} \curve(74,17, 93,17)
\put(73,16.1){$\bullet$} \put(85,13.5){$v_1$} \curve(72,19, 72,37)
\put(71.1,18){$\bullet$} \put(72.5,31){$v_2$} \curve(71,16, 56,1)
\put(70,15){$\bullet$} \put(60,2){$v_3$} \put(35,-2){Figure 3.2}
\end{picture}
\end{figure}
It is obvious that $\Gamma$ can be obtained by gluing the three half
lines emanating from the point $\xi\in N_\mathbb{R}$ in the
directions $v_1,v_2,v_3$. See Figure 3.2. These half lines are the
tropical discs which are corresponding to the three families of
Maslov index two holomorphic discs $\varphi_1,\varphi_2,\varphi_3$
respectively. We see that the above quantum relation corresponds
exactly to the equation
$$\Psi_1\star\Psi_2\star\Psi_3=q$$
in $\mathbb{C}[\Psi_1^{\pm1},\Psi_2^{\pm1}]$.

Without the assumption that $\bar{X}$ is a product of projective
spaces, the tropical interpretation will break down. This is because
for general toric Fano manifolds, the holomorphic curves which
contribute to the small quantum product may have components mapped
into the anticanonical toric divisor $D_\infty$. An example is
provided by the exceptional curve in the blowup of $\mathbb{C}P^2$
at one $T_N$-invariant point (see Example 3 in Section 4 in
\cite{Chan-Leung08}). Now the problem is that tropical geometry
cannot be used to count these holomorphic curves. In other words,
there are no tropical curves corresponding to such holomorphic
curves (cf. Rau \cite{Rau08}).\\

Now it's time to return to the main theme of this section, namely,
we can construct and apply SYZ mirror transformations to the study
of mirror symmetry for toric Fano manifolds. First we equip
$LX=X\times N$ with the symplectic structure $\pi^*(\omega_X)$,
which we denote again by $\omega_X$. Also let
$\mu_{LX}:LX\rightarrow P$ be the composition map $\mu_X\circ\pi$.
Analog to the semi-flat case, we consider the fiber product
$LX\times_P Y=P\times N\times\sqrt{-1}(T_N\times T_M)$ of the
fibrations $\mu_{LX}:LX\rightarrow P$ and $\nu_Y:Y\rightarrow P$.
Note that we have a covering map $LX\times_P Y\rightarrow X\times_P
Y$. Pulling back the universal curvature two-form
$F=\sqrt{-1}\sum_{j=1}^n dy_j\wedge du_j\in\Omega^2(X\times_P Y)$,
we get a two-form on $LX\times_P Y$, which we again denote by $F$.
We further define the \emph{holonomy function}
$\textrm{hol}:LX\times_P Y\rightarrow\mathbb{C}$ by
$$\textrm{hol}(p,v,z)=\textrm{hol}_{\nabla_y}(v)=e^{-\sqrt{-1}\langle y,v\rangle}$$
for $(p,v)\in LX, z=\exp(-x-\sqrt{-1}y)\in Y$ such that
$\mu_{X}(p)=\nu_Y(z)=x$. The SYZ mirror transformation for toric
Fano manifolds is constructed as \textit{a combination of the
semi-flat SYZ transformation $\mathcal{F}^{\textrm{sf}}$ and
fiberwise Fourier series}.
\begin{defn}
The SYZ mirror transformation
$\mathcal{F}:\Omega^*(LX)\rightarrow\Omega^*(Y)$ for $\bar{X}$ is
defined by
\begin{eqnarray*}
\mathcal{F}(\alpha) & = &
(-2\pi\sqrt{-1})^{-n}\pi_{Y,*}(\pi_{LX}^*(\alpha)\wedge e^{\sqrt{-1}F}\textrm{hol})\\
& = & (-2\pi\sqrt{-1})^{-n}\int_{N\times
T_N}\pi_{LX}^*(\alpha)\wedge e^{\sqrt{-1}F}\textrm{hol},
\end{eqnarray*}
where $\pi_{LX}:LX\times_P Y\rightarrow LX$ and $\pi_Y:LX\times_P
Y\rightarrow Y$ are the two natural projections.
\end{defn}
The basic properties of $\mathcal{F}$ are similar to those of other
Fourier-type transformations, and in particular, it satisfies the
inversion property with the \emph{inverse SYZ mirror transformation}
$\mathcal{F}^{-1}:\Omega^*(Y)\rightarrow\Omega^*(LX)$ defined by
\begin{eqnarray*}
\mathcal{F}^{-1}(\alpha) & = &
(-2\pi\sqrt{-1})^{-n}\pi_{LX,*}(\pi_Y^*(\alpha)\wedge e^{-\sqrt{-1}F}\textrm{hol}^{-1})\\
& = & (-2\pi\sqrt{-1})^{-n}\int_{T_M}\pi_Y^*(\alpha)\wedge
e^{-\sqrt{-1}F}\textrm{hol}^{-1}.
\end{eqnarray*}

In \cite{Chan-Leung08}, the SYZ mirror transformation was, for the
first time, used to study the appearance of the superpotential $W$
as quantum corrections. More precisely, we showed that
\begin{thm}[First part of Theorem
1.1 in \cite{Chan-Leung08}]\label{thm3.1} The SYZ mirror
transformation (or fiberwise Fourier series) of the function $\Psi$,
defined in terms of the counting of Maslov index two $J$-holomorphic
discs in the toric Fano manifold $\bar{X}$ with boundary in
Lagrangian torus fibers, gives the superpotential
$W:Y\rightarrow\mathbb{C}$ on the mirror manifold:
$$\mathcal{F}(\Psi)=W.$$
Furthermore, we can incorporate the symplectic structure $\omega_X$
to give the holomorphic volume form of the Landau-Ginzburg model
$(Y,W)$ in the sense that
\begin{eqnarray*}
\mathcal{F}(e^{\sqrt{-1}\omega_X+\Psi})=e^W\Omega_Y.
\end{eqnarray*}
Conversely, we have
$$\mathcal{F}^{-1}(W)=\Psi,\ \mathcal{F}^{-1}(e^W\Omega_Y)=e^{\sqrt{-1}\omega_X+\Psi}.$$
\end{thm}
\begin{nb}$\mbox{}$
\begin{enumerate}
\item[1.] We shall mention that the fact that the superpotential $W$ can be
computed in terms of the counting of Maslov index two holomorphic
discs in $\bar{X}$ with boundary in Lagrangian torus fibers was
originally due to Cho and Oh \cite{Cho-Oh03}. The key point of our
result is that there is an explicit Fourier-Mukai-type
transformation, namely, the SYZ mirror transformation $\mathcal{F}$,
that gives the superpotential $W$ by transforming an object (the
function $\Psi$) in the A-model of $\bar{X}$.
\item[2.] Apparently, the statements written here are slightly
different from those in Theorem 1.1 in \cite{Chan-Leung08}, but
realizing that $\Phi=\textrm{Exp }\Psi$, it is easy to see that they
are in fact the same statements.
\item[3.] The complex oscillatory integrals
$$\int_\Gamma e^W\Omega_Y$$
of the $n$-form $e^W\Omega_Y$ over Lefschetz thimbles $\Gamma\subset
Y$ (defined by the singularities of $W:Y\rightarrow\mathbb{C}$),
which satisfy certain Picard-Fuchs equations, play the role of
periods for Calabi-Yau manifolds. This is why we call $e^W\Omega_Y$
the holomorphic volume form of the Landau-Ginzburg model $(Y,W)$.
\end{enumerate}
\end{nb}
On the other hand, we also showed that the SYZ mirror transformation
(which, in this case, is fiberwise Fourier series)
$\mathcal{F}(\Psi_i)$ of the function $\Psi_i$ is nothing but the
monomial $e^{\lambda_i}z^{v_i}$ on $Y$, for $i=1,\ldots,d$. Since
the Jacobian ring $Jac(W)$ of the superpotential $W$ is generated by
the monomials $e^{\lambda_1}z^{v_1},\ldots,e^{\lambda_d}z^{v_d}$, by
Proposition~\ref{prop3.2}, the SYZ mirror transformation realizes a
natural isomorphism between the small quantum cohomology
$QH^*(\bar{X})$ and the Jacobian ring $Jac(W)$.
\begin{thm}[Second part of Theorem
1.1 in \cite{Chan-Leung08}]\label{thm3.2} The SYZ mirror
transformation $\mathcal{F}$ induces a natural isomorphism of
$\mathbb{C}$-algebras
\begin{eqnarray*}
\mathcal{F}:QH^*(\bar{X})\overset{\cong}{\longrightarrow} Jac(W),
\end{eqnarray*}
which takes the quantum product, now realized as a convolution
product, to the ordinary product of Laurent polynomials, provided
that $\bar{X}$ is a product of projective spaces.
\end{thm}
In the example of $\bar{X}=\mathbb{C}P^2$, the superpotential is the
Laurent polynomial $W(z_1,z_2)=z_1+z_2+\frac{q}{z_1z_2}$ on
$Y=(\mathbb{C}^*)^2$, where $q=e^{-t}$. Its logarithmic partial
derivatives are given by
$$\partial_1W=z_1-\frac{q}{z_1z_2},\ \partial_2W=z_2-\frac{q}{z_1z_2},$$
so that the Jacobian ring is given by
\begin{eqnarray*}
Jac(W) & = & \mathbb{C}[z_1^{\pm1},z_2^{\pm1}]\big/\big\langle
z_1-\frac{q}{z_1z_2},z_2-\frac{q}{z_1z_2}\big\rangle\\
& = & \mathbb{C}[Z_1,Z_2,Z_3]\big/\big\langle
Z_1-Z_3,Z_2-Z_3,Z_1Z_2Z_3-q\big\rangle,
\end{eqnarray*}
where the monomials $Z_1=z_1$, $Z_2=z_2$ and $Z_3=\frac{q}{z_1z_2}$
are the SYZ mirror transformations (i.e. fiberwise Fourier series)
of the functions $\Psi_1$, $\Psi_2$ and $\Psi_3$ respectively.
\begin{nb}$\mbox{}$
\begin{enumerate}
\item[1.] In \cite{CCIT06}, Coates, Corti, Iritani
and Tseng formulated the mirror symmetry conjecture for toric
manifolds (and orbifolds) as an isomorphism of graded
$\frac{\infty}{2}$VHS between the A-model $\frac{\infty}{2}$VHS
associated to a toric manifold and the B-model $\frac{\infty}{2}$VHS
associated to the mirror Landau-Ginzburg model (see also Iritani
\cite{Iritani07}). It is desirable to have this isomorphism, which
contains more information than the isomorphism in the above theorem,
realized by SYZ mirror transformations.
\item[2.]
In \cite{FOOO08a} (and also \cite{FOOO08b}), Fukaya-Oh-Ohta-Ono
applied the machinery developed in \cite{FOOO00} to the case of
toric manifolds. They considered Floer cohomology with coefficients
in the Novikov ring, instead of $\mathbb{C}$ used here and in
Auroux's paper \cite{Auroux07}. They have results on the
superpotential even in the non-Fano toric case. The isomorphism
$QH^*(\bar{X})\cong Jac(W)$ (over the Novikov ring) was also
discussed and proved in their work (Theorem 1.9 in \cite{FOOO08a}).
Their proof is combinatorial, using Batyrev's presentation of the
small quantum cohomology ring for toric Fano manifolds, the validity
of which in turn relies on Givental's mirror theorem. They claimed
that a more conceptual and geometric proof for toric, not
necessarily Fano, manifolds will appear in a sequel to their paper.
\end{enumerate}
\end{nb}

\subsection{Transformation of branes}

This subsection is an attempt to understand the correspondence
between A-branes of the toric Fano manifold $\bar{X}$ and B-branes
of the mirror Landau-Ginzburg model $(Y,W)$ via SYZ mirror
transformations.

We will deal with the simplest case of the correspondence. So let
$L_x=\mu_X^{-1}(x)$ be the Lagrangian torus fiber of $\bar{X}$ over
a point $x\in P$. We equip $L_x$ with a flat $U(1)$-bundle
$\mathbb{L}_y=(L_x\times\mathbb{C},\nabla_y)$, where $\nabla_y$ is
the flat $U(1)$-connection corresponding to $y\in(L_x)^\vee$. The
mirror of the A-brane $(L_x,\mathbb{L}_y)$ is given, according to
SYZ, by the B-brane $(z=\exp(-x-\sqrt{-1}y)\in Y,\mathcal{O}_z)$. In
other words, the correspondence on the level of objects is the same
as in the semi-flat Calabi-Yau case. Quantum corrections will emerge
and make a difference when we consider their endomorphisms.

According to Hori (see~\cite{HKKPTVVZ03}, Chapter 39), the
endomorphism algebra $\textrm{End}(z,\mathcal{O}_z)$ of the B-brane
$(z,\mathcal{O}_z)$, as a $\mathbb{C}$-vector space, is given by the
cohomology of the complex
$$({\bigwedge}^{\!*}T_zY,\delta=\iota_{\partial W(z)}),$$
where $\iota_{\partial W(z)}$ is contraction with the vector
$\partial W(z)=\sum_{j=1}^n\partial_jW(z)(\partial_j)_z$ and here
again $\partial_j$ denotes $z_j\frac{\partial}{\partial z_j}$. The
following elementary proposition shows that the introduction of the
superpotential $W$ "localizes" the category B-branes to the critical
points of $W$.
\begin{prop}
The endomorphism $\textrm{End}(z,\mathcal{O}_z)$ is nontrivial if
and only if $z\in Y$ is a critical point of the superpotential
$W:Y\rightarrow\mathbb{C}$, and in which case,
$\textrm{End}(z,\mathcal{O}_z)$ is isomorphic to
${\bigwedge}^{\!*}T_zY$ as $\mathbb{C}$-vector spaces.
\end{prop}
On the other hand, the endomorphism algebra of the A-brane
$(L_x,\mathbb{L}_y)$ in the (derived) Fukaya category is given by
the \emph{Floer cohomology ring} $HF(L_x,\mathbb{L}_y)$,\footnote{We
use $\mathbb{C}$ as the coefficient ring, instead of the Novikov
ring.} which in turn, as a $\mathbb{C}$-vector space, is given by
the cohomology of the Floer complex
$$(C^*(L_x,\mathbb{C}),\delta=m_1)$$
where $m_1=m_1(L_x,\mathbb{L}_y)$ denotes the Floer differential. In
\cite{Cho-Oh03}, \cite{Cho04}, Cho and Oh explicitly computed the
Floer differential $m_1$. Recall that $H^1(L_x,\mathbb{C})$, viewed
as the space of infinitestimal deformations of the pair
$(L_x,\mathbb{L}_y)$, is canonically isomorphic to $T_zY$. Let
$C_1,\ldots,C_n$ be the basis of $H^1(L_x,\mathbb{C})$ corresponding
to $(\partial_1)_z,\ldots,(\partial_n)_z$. Then the results of Cho
and Oh stated that
$m_{1,\beta_i}(C_j)=C_j\cdot\partial\beta_i=v_i^j$ and
\begin{eqnarray*}
m_1(C_j) & = & \sum_{i=1}^dm_{1,\beta_i}(C_j)
\exp(-\frac{1}{2\pi}\int_{\beta_i}\omega_X)\textrm{hol}_{\nabla_y}(\partial\beta_i)\\
& = & \sum_{i=1}^dv_i^jz^{v_i}=\partial_jW(z).
\end{eqnarray*}
This shows that $m_1=\iota_{\partial W(z)}$ on
$H^1(L_x,\mathbb{C})=T_zY$, and $m_1=0$ on $H^1(L_x,\mathbb{C})$ if
and only if $z$ is a critical point of $W$. The following result
proved by Cho-Oh in \cite{Cho-Oh03} is parallel to the above
proposition.
\begin{thm}[Cho-Oh \cite{Cho-Oh03}]
The Floer cohomology $HF(L_x,\mathbb{L}_y)$ is nontrivial and
isomorphic to $H^*(L_x,\mathbb{C})$ if and only if $m_1=0$ on
$H^1(L_x,\mathbb{C})$.
\end{thm}
We conclude that
\begin{thm}\label{thm3.4}
The Floer cohomology $HF(L_x,\mathbb{L}_y)$ of the A-brane
$(L_x,\mathbb{L}_y)$ is isomorphic to the endomorphism algebra
$\textrm{End}(z,\mathcal{O}_z)$ of the mirror B-brane
$(z,\mathcal{O}_z)$ as $\mathbb{C}$-vector spaces.
\end{thm}
It is intriguing to see whether this isomorphism can be realized by
explicit SYZ mirror transformations.
\begin{nb}
In \cite{Cho04}, Cho proved that the Floer cohomology ring
$HF(L_x,\mathbb{L}_y)$, equipped with the product structure given by
$m_2=m_2(L_x,\mathbb{L}_y)$, is a Clifford algebra generated by
$H^1(L_x,\mathbb{C})$ with the bilinear form given by the Hessian of
$W$: $Q(C_j,C_k)=\partial_j\partial_kW(z)$. This implies that the
isomorphism in Theorem~\ref{thm3.4} is in fact an isomorphism of
$\mathbb{C}$-algebras. This confirms a prediction by Physicists. See
the paper of Cho \cite{Cho04} for details.
\end{nb}

\section{Further questions}

The results described in this article represent the first step in
our program which is aimed at exploring mirror symmetry via SYZ
mirror transformations. In particular, they showed that these
transformations can be applied successfully to explain the mirror
symmetry for toric Fano manifolds, a case where quantum corrections
do exist. However, we shall emphasize that the quantum corrections
in the toric Fano case, which are due to the anticanonical toric
divisor, are much simpler than those in the general case
(Gross-Siebert \cite{GS07}, Auroux \cite{Auroux07}), where quantum
corrections may arise due to contributions from the proper singular
Lagrangian fibers of the Lagrangian torus fibrations and complicated
wall-crossing phenomena start to interfere. In terms of affine
geometry, this means that the bases of the Lagranigan torus
fibrations in the toric case are affine manifolds with boundary but
without singularities, while in the general case, the bases are
affine manifolds with \textit{both} boundary and singularities (and
in the semi-flat case, the bases are affine manifolds without
boundary and singularities). Certainly much more work remains to be
done in the future. In this final section, we will comment on
several possible future research directions. The discussion is going
to be rather speculative.

\subsection{Toric Fano manifolds}

We have seen that the simplest correspondence between A-branes on a
toric Fano manifold $\bar{X}$ and B-branes on the mirror
Landau-Ginzburg model $(Y,W)$, namely
$$(L_x,\mathbb{L}_y)\longleftrightarrow(z,\mathcal{O}_z),$$
is compatible with the SYZ philosophy. It is desirable to see how
other A-branes on $X$ are transformed to the corresponding mirror
B-branes on $(Y,W)$. An interesting and important example would be
the Lagrangian submanifold $\mathbb{R}P^n\subset\mathbb{C}P^n$ for
odd $n$, which can be viewed as a multi-section of the moment map of
$\mathbb{C}P^n$. Employing the SYZ approach, the mirror B-brane is
expected to be a trivial rank-$2^n$ holomorphic vector bundle over
$Y$, equipped with some additional information related to $W$. A
possible choice of this additional information would be a
\emph{matrix factorization} of $W$; currently, it is widely believed
that the category of B-branes on $(Y,W)$ is given by the category of
matrix factorizations of $W$. This was first proposed by Kontsevich,
see Orlov \cite{Orlov03} for details. The relation between these
matrix factorizations and the computation of Floer cohomology will
be the key to a complete understanding of the correspondences of
branes.\\

On the other hand, we have not even touched the correspondence
between B-branes on $\bar{X}$ and A-branes on $(Y,W)$. As we
mentioned in the introduction, the results of Seidel
\cite{Seidel00}, Ueda \cite{Ueda04}, Auroux-Katzarkov-Orlov
\cite{AKO04}, \cite{AKO05} and Abouzaid \cite{Abouzaid05},
\cite{Abouzaid06} have provided substantial evidences for this half
of the Homological Mirror Symmetry Conjecture. In particular,
Abouzaid \cite{Abouzaid06} made use of an idea originated from the
SYZ conjecture, namely, the mirror of a Lagrangian section should be
a holomorphic line bundle. His results also showed that the
correspondence is in line with the SYZ picture. Recently, Fang
\cite{Fang08} and Fang-Liu-Treumann-Zaslow \cite{FLTZ08} proved a
version of Homological Mirror Symmetry for toric manifolds by
explicitly using T-duality. It is an interesting question whether
one can construct an explicit SYZ mirror transformation to realize
the correspondence between B-branes on $\bar{X}$ and A-branes on
$(Y,W)$.

\subsection{Toric non-Fano or non-toric Fano manifolds}

As in the case of toric Fano manifolds, non-toric Fano manifolds
such as Grassmannians and flag manifolds admit natural Lagrangian
torus fibrations, provided by Gelfand-Cetlin integrable systems
(see, for example, Guillemin-Sternberg \cite{GS83}), which are
convenient for applying SYZ mirror transformations. While mirror
symmetry for these manifolds has been studied for some time by
Givental \cite{Givental96} and others, new tools and new ideas are
needed if we want to apply SYZ mirror transformations to these
examples. The recent works of Nishinou-Nohara-Ueda \cite{NNU08a},
\cite{NNU08b} have shed some light on this case. In particular, they
obtained a classification the holomorphic discs in flag manifolds
with boundary in Lagrangian torus fibers, which should be very
useful in the constructions of SYZ mirror transformations.

On the other hand, the mirror symmetry for toric non-Fano manifolds
is also not well understood too. As can be seen from the works of
Givental \cite{Givental97}, the mirror map between the complexified
K\"{a}hler and complex moduli spaces in this case is a nontrivial
\emph{coordinate change}, instead of an identity map as in the toric
Fano case. In Auroux \cite{Auroux07}, nontrivial coordinate changes
and wall-crossing phenomena were also observed in constructing the
superpotentials for the mirrors of non-toric examples. Hence, the
definitions of the SYZ mirror transformations may have to be
adjusted to incorporate the nontrivial mirror map and also
wall-crossing phenomena. For this, we will have to make the
construction of SYZ mirror transformations local. A very preliminary
attempt to this is made in Section 5 in \cite{Chan-Leung08}.

\subsection{Calabi-Yau manifolds}

The ultimate goal of our program is no doubt to apply SYZ mirror
transformations to get a better understanding of the mirror symmetry
for Calabi-Yau manifolds and the SYZ Conjecture. Works of Fukaya
\cite{Fukaya02}, Kontsevich-Soibelman \cite{KS04} and Gross-Siebert
\cite{GS07} have laid an important foundation for understanding the
SYZ framework for both Calabi-Yau and non-Calabi-Yau manifolds. In
view of the fact that toric varieties have played an important role
in the constructions of Gross and Siebert, it would be nice if we
can incorporate our methods with their new techniques to study SYZ
mirror transformations for Calabi-Yau manifolds; and hopefully, this
would let us reveal geometrically the secret of mirror symmetry.



\begin{thebibliography}{99}


\bibitem{Abouzaid05}
M. Abouzaid, Homogeneous coordinate rings and mirror symmetry for
toric varieties, Geom. Topol., {\bf 10} (2006), 1097--1157
(math.SG/0511644).

\bibitem{Abouzaid06}
\underline{\qquad\quad}, Morse homology, tropical geometry, and
homological mirror symmetry for toric varieties, preprint, 2006
(math/0610004).

\bibitem{Abreu00}
M. Abreu, K\"{a}hler geometry of toric manifolds in symplectic
coordinates, Symplectic and contact topology: interactions and
perspectives (Toronto, ON/Montreal, QC, 2001), 1--24, Fields Inst.
Commun., {\bf 35}, Amer. Math. Soc., Providence, RI, 2003
(math.DG/0004122).

\bibitem{Auroux07}
D. Auroux, Mirror symmetry and T-duality in the complement of an
anticanonical divisor, J. Gokova Geom. Topol. GGT, {\bf 1} (2007),
51--91 (arXiv:0706.3207).

\bibitem{AKO04}
D. Auroux, L. Katzarkov and D. Orlov, Mirror symmetry for weighted
projective planes and their noncommutative deformations, Ann. of
Math. (2), {\bf 167} (2008), no. 3, 867--943 (math.AG/0404281).

\bibitem{AKO05}
\underline{\qquad\quad}, Mirror symmetry for del Pezzo surfaces:
vanishing cycles and coherent sheaves, Invent. Math., {\bf 166}
(2006), no. 3, 537--582 (math.AG/0506166).

\bibitem{Chan-Leung08}
K.-W. Chan and N.-C. Leung, Mirror symmetry for toric Fano manifolds
via SYZ transformations, preprint, 2008 (arXiv:0801.2830).

\bibitem{Cho04}
C.-H. Cho, Products of Floer cohomology of torus fibers in toric
Fano manifolds, Comm. Math. Phys., {\bf 260} (2005), no. 3, 613--640
(math.SG/0412414).

\bibitem{Cho-Oh03}
C.-H. Cho and Y.-G. Oh, Floer cohomology and disc instantons of
Lagrangian torus fibers in Fano toric manifolds, Asian J. Math.,
{\bf 10} (2006), no. 4, 773--814 (math.SG/0308225).

\bibitem{CCIT06}
T. Coates, A. Corti, H. Iritani and H.-H. Tseng, Wall-crossings in
toric Gromov-Witten theory I: crepant examples, preprint, 2006
(math.AG/0611550).

\bibitem{Fang08}
B. Fang, Homological mirror symmetry is T-duality for
$\mathbb{P}^n$, preprint, 2008 (arXiv:0804.0646).

\bibitem{FLTZ08}
B. Fang, C.-C. M. Liu, D. Treumann and E. Zaslow, T-duality and
equivariant homological mirror symmetry for toric varieties,
preprint, 2008 (arXiv:0811.1228).

\bibitem{Fukaya02}
K. Fukaya, Multivalued Morse theory, asymptotic analysis and mirror
symmetry, Graphs and patterns in mathematics and theoretical
physics, 205--278, Proc. Sympos. Pure Math., {\bf 73}, Amer. Math.
Soc., Providence, RI, 2005.

\bibitem{FOOO00}
K. Fukaya, Y.-G. Oh, H. Ohta and K. Ono, Lagrangian intersection
Floer theory: Anomaly and obstruction, preprint, 2000.

\bibitem{FOOO08a}
\underline{\qquad\quad}, Lagrangian Floer theory on compact toric
manifolds I, preprint, 2008 (arXiv:0802.1703).

\bibitem{FOOO08b}
\underline{\qquad\quad}, Lagrangian Floer theory on compact toric
manifolds II: Bulk deformations, preprint, 2008 (arXiv:0810.5654).

\bibitem{Givental94}
A. Givental, Homological geometry and mirror symmetry, Proceedings
of the International Congress of Mathematicians, Vol. {\bf 1, 2}
(Z\"{u}rich, 1994), 472--480, Birkh\"{a}user, Basel, 1995.

\bibitem{Givental96}
\underline{\qquad\quad}, Stationary phase integrals, quantum Toda
lattices, flag manifolds and the mirror conjecture, Topics in
singularity theory, 103--115, Amer. Math. Soc. Transl. Ser. 2, {\bf
180}, Amer. Math. Soc., Providence, RI, 1997 (alg-geom/9612001).

\bibitem{Givental97}
\underline{\qquad\quad}, A mirror theorem for toric complete
intersections, Topological field theory, primitive forms and related
topics (Kyoto, 1996), 141--175, Progr. Math., {\bf 160},
Birkh\"{a}user Boston, Boston, MA, 1998 (alg-geom/9701016).

\bibitem{Gross09}
M. Gross, Mirror symmetry for $\mathbb{P}^2$ and tropical geometry,
preprint, 2009 (arXiv:0903.1378).

\bibitem{GS07}
M. Gross and B. Siebert, From real affine geometry to complex
geometry, preprint, 2007 (math.AG/0703822).

\bibitem{Guillemin94}
V. Guillemin, Moment maps and combinatorial invariants of
Hamiltonian $T^n$-spaces, Progress in Mathematics, {\bf 122}.
Birkha\"{u}ser Boston, Inc., Boston, MA, 1994.

\bibitem{GS83}
V. Guillemin and S. Sternberg, \textit{The Ge\'{l}fand-Cetlin system
and quantization of the complex flag manifolds}, J. Funct. Anal. 52
(1983), no. 1, 106--128.

\bibitem{Hitchin97}
N. Hitchin, The moduli space of special Lagrangian submanifolds,
Ann. Scuola Norm. Sup. Pisa Cl. Sci. (4) {\bf 25} (1997), no. 3-4,
503--515 (dg-ga/9711002).

\bibitem{HIV00}
K. Hori, A. Iqbal and C. Vafa, D-branes and mirror symmetry,
preprint, 2000 (hep-th/0005247).

\bibitem{HKKPTVVZ03}
K. Hori, S. Katz, A. Klemm, R. Pandharipande, R. Thomas, C. Vafa, R.
Vakil, and E. Zaslow, Mirror symmetry, Clay Mathematics Monographs,
{\bf 1}, American Mathematical Society, Providence, RI; Clay
Mathematics Institute, Cambridge, MA, 2003.

\bibitem{Hori-Vafa00}
K. Hori and C. Vafa, Mirror symmetry, preprint, 2000
(hep-th/0002222).

\bibitem{Iritani07}
H. Iritani, Real and integral structures in quantum cohomology I:
toric orbifolds, preprint, 2007 (arXiv:0712.2204).

\bibitem{Kontsevich98}
M. Kontsevich, Lectures at ENS, Paris, Spring 1998, notes taken by
J. Bellaiche, J.-F. Dat, I. Marin, G. Racinet and H.
Randriambololona.

\bibitem{KS04}
M. Kontsevich and Y. Soibelman, Affine structures and
non-Archimedean analytic spaces, The unity of mathematics, 321--385,
Progr. Math., {\bf 244}, Birkh\"{a}user Boston, Boston, MA, 2006
(math.AG/0406564).

\bibitem{Leung00}
N.-C. Leung, Mirror symmetry without corrections, Comm. Anal. Geom.
{\bf 13} (2005), no. 2, 287--331 (math.DG/0009235).

\bibitem{LYZ00}
N.-C. Leung, S.-T. Yau and E. Zaslow, From special Lagrangian to
Hermitian-Yang-Mills via Fourier-Mukai transform, Adv. Theor. Math.
Phys., {\bf 4} (2000), no. 6, 1319--1341 (math.DG/0005118).

\bibitem{Mikhalkin03}
G. Mikhalkin, \textit{Enumerative tropical algebraic geometry in
$\mathbb{R}^2$}, J. Amer. Math. Soc. 18 (2005), no. 2, 313--377
(math.AG/0312530).

\bibitem{NNU08a}
T. Nishinou, Y. Nohara and K. Ueda, Toric degenerations of
Gelfand-Cetlin systems and potential functions, preprint, 2008
(arXiv:0810.3470).

\bibitem{NNU08b}
\underline{\qquad\quad}, Potential functions via toric
degenerations, preprint, 2008 (arXiv:0812.0066).

\bibitem{NS04}
T. Nishinou and B. Siebert, \textit{Toric degenerations of toric
varieties and tropical curves}, Duke Math. J. 135 (2006), no. 1,
1--51 (math.AG/0409060).

\bibitem{Orlov03}
D. Orlov, Triangulated categories of singularities and D-branes in
Landau-Ginzburg models, Proc. Steklov Inst. Math. 2004, no. 3, {\bf
246}, 227--248 (math.AG/0302304).

\bibitem{Rau08}
J. Rau, Intersections on tropical moduli spaces, preprint, 2008
(arXiv:0812.3678).

\bibitem{Seidel00}
P. Seidel, More about vanishing cycles and mutation, Symplectic
geometry and mirror symmetry (Seoul, 2000), 429--465, World Sci.
Publ., River Edge, NJ, 2001 (math.SG/0010032).

\bibitem{SYZ96}
A. Strominger, S.-T. Yau and E. Zaslow, Mirror symmetry is
T-duality, Nuclear Phys. B, {\bf 479} (1996), no. 1-2, 243--259
(hep-th/9606040).

\bibitem{Ueda04}
K. Ueda, Homological mirror symmetry for toric del Pezzo surfaces,
Comm. Math. Phys., {\bf 264} (2006), no. 1, 71--85
(math.AG/0411654).

\end{thebibliography}
\end{document}